\documentclass[reqno]{amsart}
\usepackage{bbold}
\usepackage{color}
\usepackage{hyperref}
\usepackage{color}
\usepackage[latin1]{inputenc}
\usepackage{epsfig}
\usepackage{amssymb}
\usepackage{color}
\usepackage{dsfont}
\usepackage{amsthm}
\usepackage{amsmath}
\usepackage{amsfonts}
\usepackage{amssymb}
\usepackage{cancel}
\usepackage{nicefrac}
\usepackage{graphicx}
\usepackage{epsfig}
\usepackage{floatrow}
\usepackage{lmodern}
\usepackage{ae}

\newtheorem{lem}{Lemma}[section]

\newtheorem{thm}{Theorem}[section]
\newtheorem{prop}{Proposition}[section]
\newtheorem{defi}{Definition}[section]

\theoremstyle{remark}

\newfont{\sBlackboard}{msbm10 scaled 900}

\newcommand{\mylabel}[1]{\label{#1}
            \ifx\undefined\stillediting
            \else \fbox{$#1$}\fi }
\newcommand{\BE}{\begin{equation}}

\newcommand{\EEQ}{\end{equation}}
\newcommand{\rfb}[1]{\mbox{\rm
   \eqref{#1}}\ifx\undefined\stillediting\else:\fbox{$#1$}\fi}

\newfont{\Blackboard}{msbm10 scaled 1200}

\newfont{\roma}{cmr10 scaled 1200}

\def\n{|\kern -.05cm{|}\kern -.05cm{|}}

\def\N{\rm I\hskip -2pt N} 
\def\Z{{\bf Z}} 

\newcommand{\mm}    {{\hbox{\hskip 0.5pt}}}

\newcommand{\bluff} {{\hbox{\raise 15pt \hbox{\mm}}}}
scaled\magstep2

\makeatletter
\def\section{\@startsection {section}{1}{\z@}{-3.5ex plus -1ex minus
    -.2ex}{2.3ex plus .2ex}{\large\bf}}
\makeatother
\def\be{\begin{equation}}
\def\ee{\end{equation}}

\begin{document}
\thispagestyle{empty}
\title[Qualitative properties for  $1-D$ impulsive wave equation]{Qualitative properties for the $1-D$ impulsive wave equation: controllability and observability}
\author{Akram BEN  AISSA$^{*}$}
\author{Walid Zouhair$^{**}$}

\begin{abstract} 
In this paper, we establish some important results for the impulsive wave equation. We begin by proving the existence of a solution. Then, we study the impulse approximate controllability where the control function acts on a subdomain $ \omega $ and at one instant of time $ \tau \in (0, T)$. Afterward, we study impulse observability.
 \end{abstract}


\keywords{impulsive wave equation, impulse control, impulse observability, impulse approximate controllability. \\
*: UR Analysis and Control of PDE's, UR 13ES64, Higher Institute of transport and Logistics of Sousse, University of Sousse,  Tunisia\\
akram.benaissa@fsm.rnu.tn\\
**: Cadi Ayyad University, Faculty of Sciences Semlalia, LMDP, UMMISCO (IRD-UPMC), B.P. 2390, Marrakesh, Morocco,\\ walid.zouhair.fssm@gmail.com}

\maketitle

\section{Introduction}
Impulsive dynamic systems are a type of hybrid systems for which the trajectory admits discontinuities at certain instants due to sudden jumps of the state called pulses ( see more in \cite{BMJ}). Dynamical behavior of many systems in real life can be characterized by abrupt changes that appear suddenly, such as heartbeats, drug flows, the value of stocks, impulse vaccination, and bonds on the stock market. This class of hybrid systems is  presented as follows:

\begin{equation}\label{S1}
\left\{
\begin{array}{lll}
\xi^{\prime}(t) &=& f(t,\xi(t)),\quad  t \in (0,T],\,\, t \neq t_{k},\\[3mm]
\bigtriangleup \xi(t_{k})&=& \nu_{k}(\xi(t_{k})),\quad k \in \theta_{m}^{n},
\end{array}
\right.
\end{equation}
where $\theta_{m}^{n} = \lbrace m,m+1,...,n \rbrace ,$ \, $\bigtriangleup \xi(t_{k}) = \xi (t_{k}^{+})- \xi (t_{k}^{-}),$ with $\xi (t_{k}^{+})$ $(\text{respectively}\, \xi (t_{k}^{-}))$ denotes the limit to the right (respectively to the left) of $t$ and $\xi^{\prime}(t)= f(t,\xi(t)),$  is a differential equation.\\
\begin{center}
\includegraphics[scale=0.5]{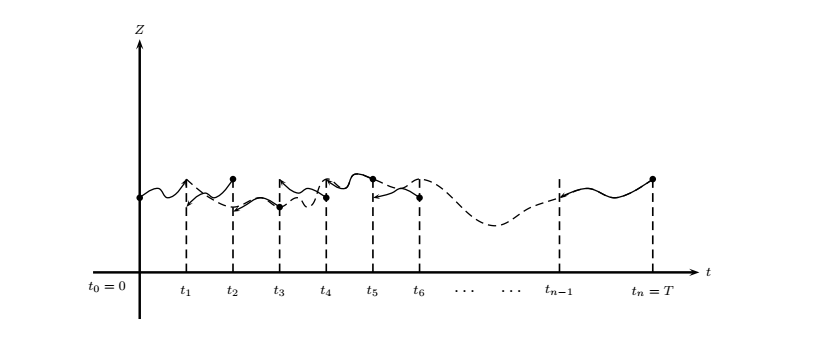}
\textbf{Explanatory diagram for the behavior of the solution }\par\medskip
\end{center}

The theory of impulsive differential equations was initiated by V. D. Mil'man and A. Mishkis in 1960 \cite{MVDM}. Afterward, many scientists contributed to the enrichment of this theory, they launched different studies on this discipline and a large number of results were established, we invite the interested reader to see \cite{liu,farzana,yao}. \\

Unlike the interior controllability for impulsive systems that have been extensively studied in the literature, see for instance \cite{CDJUHLOC,ACCDGHL,GL,CDGHL,LHWZME} and the references therein, the problem of controllability with impulse controls has attracted less attention and not as many works are available in this area. We mention A. Khapalov \cite{KAY} who proved the exact controllability of a class of second-order hyperbolic boundary problems with impulse controls using the Huygens' principle. S. Qin and G. Wang \cite{SQG} presented a necessary and sufficient condition for the approximate controllability in finite dimension, which is exactly Kalman's controllability rank condition. Recently K. D. Phung, G. Wang, and Y. Xu \cite{KDP} considered the following impulsive functional differential systems 
\begin{equation}
\left\{\begin{array}{ll}
\partial_{t} \psi-\Delta \psi=0, & \text { in } \Omega \times(0, T) \backslash\{\tau\} \\
\psi=0, & \text { on } \partial \Omega \times(0, T) \\
\psi(\cdot, 0)=\psi^{0}, & \text { in } \Omega \\
\psi(\cdot, \tau)=\psi\left(\cdot, \tau^{-}\right)+\mathbb{1}_{\omega} h, & \text { in } \Omega.
\end{array}\right.
\end{equation}
They proved that the above system is null approximate impulse controllable at any time $T > 0$.

Motivated by the above works, we study in this paper the following  $1-D$ impulsive wave equation:
\begin{equation}\label{s2}
\left\{
\begin{array}{lll}
\frac{\partial^{2}\kappa(x,t)}{\partial{t}^{2}}-\frac{\partial^{2}\kappa(x,t)}{\partial{x}^{2}} = 0, &\text{ in }\,\,\,\, \Omega\, \times(0, T) \backslash\left\{t_{k}\right\}_{k \in \theta_{1}^{n}},& \\
\kappa(x,t_{k}^{+}) - \kappa(x,t_{k}^{-}) = 0, &\text{ in }\,\,\,\, \Omega, \,\,\,\,k \in \theta_{1}^{n},&\\
\frac{\partial \kappa(x,t_{k}^{+})}{\partial{t}} -\frac{\partial \kappa(x,t_{k}^{-})}{\partial{t}} = \Upsilon_{k}(x), &\text{ in }\,\,\,\, \Omega, \,\,\,\,k \in \theta_{1}^{n},&\\
\kappa(0,t) = \kappa(1,t)= 0, & \text{ on }\,\,\,\,(0, T),&\\
\kappa(x,0) = \kappa^{0}(x) , \frac{\partial \kappa(x,0)}{\partial{t}}=\kappa^{1}(x),  &\text{ in }\,\,\,\, \Omega,& \;\;\;\;
\end{array}
\right.
\end{equation}
where $\Omega = (0,1),$  $\displaystyle\lbrace t_{k}\rbrace_{k \in \theta_{1}^{n}}\,$ are instantaneous pulses such that $t_{0}=0<t_{1}<\dots<t_{n}<T$ and $\Upsilon_{k}$ is an impulse control that satisfies a suitable assumptions, $(\kappa^{0},\kappa^{1})$ are the initial conditions in an appropriate Hilbert space $\mathcal{H}.$\\

The rest of the paper is organized as follows: In section 2, we briefly present the problem formulation and the well-posedness result. In section 3 we discuss the impulse approximate controllability and we present some related definitions and propositions. The last section is devoted to studying impulse observability.

\section{Well-posedness}
In order to use the theory of semigroup to establish the well-posedness of problem \eqref{s2}, we start by putting the previous equation to the first order in time $\psi = \left( 
\kappa,\frac{\partial \kappa}{\partial{t}}\right)^\top$. Then, the problem \eqref{s2} turns into a Cauchy problem as 
\begin{equation}\label{s3}
\left\{
\begin{array}{lll}
\frac{\partial \Psi(.,t)}{\partial{t}} =\mathcal{A}\Psi(.,t) \quad & \text{ in }\,\,\,\,(0, T) \backslash\left\{t_{k}\right\}_{k \in \theta_{1}^{n}},&\\
\Psi(t_{k}^{+})=\Psi(t_{k}^{-})+\Upsilon_{k}(.)B \quad & k \in \theta_{1}^{n},&\\
\Psi(0) = \Psi_{0}.  \\
\end{array}
\right.
\end{equation}
First we introduce the Hilbert space as follows:
$$\mathcal{H}= H^{1}_{0}(\Omega)\times L^{2}(\Omega),$$  endowed with the following inner product
\begin{equation}\label{classicinnerproduct}
\left( \left(u, v \right)^\top,\,\left(   \tilde{u},\tilde{v}\right)^\top\right)_{\mathcal{H}}=
 \int_\Omega \frac{\partial u}{\partial x} \frac{\overline{\tilde{u}}}{\partial x}dx+\int_\Omega v\overline{\tilde{v}}dx.
 \end{equation}
Then, we define $B$ and the operator $\mathcal{A}$  on $\mathcal{H}$ by
 $$\mathcal{D}(\mathcal{A}) = H^{2}(\Omega)\cap H^{1}_{0}(\Omega)\times H^{1}_{0}(\Omega)$$
 and
 \begin{eqnarray}
 \mathcal{A}\left( 
u,v \right)^\top&=&\left(
v, \frac{\partial^2 u}{\partial x^2}
 \right)^\top, \forall (u,\,v)^\top\in D(\mathcal{A})\label{fi}\\
 B &=&\left(0,1 \right)^\top. \nonumber
 \end{eqnarray}
One can easily check that $\mathcal{A}$ is an unbounded skew-adjoint operator on $\mathcal{H}$, then it generates a group of isometries  $(S(t))_{t\in \mathbb{R}}\,$( Stone's theorem\cite{pazy} ).\\

For further purposes and with the notations $\mathcal{I}:=[0,T],$ and $\mathcal{I}^{\prime}:=[0,T]\backslash\left\{t_{k}\right\}_{k \in \theta_{1}^{n}},$ we set the following Banach space which is natural framework space for evolution equations with pulses:
\begin{equation}
\displaystyle \mathcal{P C}(\mathcal{I} ; \mathcal{H})=\left \{\begin{array}{ll}
y:\mathcal{I} \rightarrow \mathcal{H}: y \in C(\mathcal{I}^{\prime};\mathcal{H}):\forall k\quad y\left(t_{k}^{-}\right), y\left(t_{k}^{+}\right)\text {exists}\\ \mbox{ and }y\left(t_{k}\right)= y\left(t_{k}^{-}\right)
\end{array}\right\},
\end{equation}
endowed with the norm
\begin{equation} 
 \|y\| = \sup_{t \in [0,T]} \|y(t)\|.
 \end{equation}
We define the classical solution for the impulsive wave equation \eqref{s3} as follows:

\begin{defi}
Classical solution $\Psi$ for \eqref{s3}is an absolutely continuous piecewise continuous mapping with discontinuities at points $t = t_{k}$ which, for almost all $t$ satisfies the system \eqref{s3}, and for $t=t_{k}$ satisfies the jump conditions. In other words, a classical solution for (\ref{s3}) is a function 
\begin{equation*}
\left\{
\begin{array}{lll}
\Psi \in \mathcal{P} \mathcal{C}([0, T] ; \mathcal{H}) \,\cap\, C^{1}\left((0, T) \backslash\left\{t_{k}\right\}_{k \in \theta_{1}^{n}}, \mathcal{H}\right),\\ 
\Psi(t) \in D(A), \text { for } t \in(0, T) \backslash\left\{t_{k}\right\}_{k \in \theta_{1}^{n}}
\end{array}\right.
\end{equation*}
such that $\Psi$ satisfies \eqref{s3} in $[0,T).$
\end{defi}

\begin{prop}[Well-posedness]
Assume that\, $\Psi_{0}= (\kappa^{0},\kappa^{1}) \in \mathcal{D}(\mathcal{A})\, ,$ and $\Upsilon_{k} \in H^{1}_{0}(\Omega)$ $ k \in \theta_{1}^{n}.$ Then, the impulsive system \eqref{s3}  has a unique classical solution $\Psi$ which, for $t \in [0,T),$ satisfies
$$\Psi(t) = S(t) \Psi_{0}+\displaystyle{\sum_{0<t_{k}<t}} S(t-t_{k})\Upsilon_{k} B .$$

\end{prop}
\begin{proof}
First,  in the interval $J_1 = [0, t_{1}),$ we consider the equation
\begin{equation*}
\left\{
\begin{array}{lll}
\frac{\partial \Psi(.,t)}{\partial{t}} = \mathcal{A}\Psi(.,t), \quad \, 0<t<t_{1},\\[2mm]
\Psi(0) = \Psi_{0}.  \\
\end{array}
\right.
\end{equation*}
\\
Note that the classical solutions for evolution equations without pulses are defined in an obvious way, see Pazy \cite{pazy}. Then it is clear that the unique classical solution of this system is given by\\
$$ \Psi_{1}(t) = S(t)\Psi_{0}.$$
Next, we define
$$ \Psi_{1}(t_{1}) = S(t_{1})\Psi_{0},$$ 
furthermore, we can check that $\Psi_{1}(.) $ is left continuous at $t_{1},$ and $\Psi_{1}(t_{1}) \in \mathcal{D}(\mathcal{A}). $\\
On the other hand in $J_2 = [t_{1}, t_{2})$, we  consider the following equation
\begin{equation*}
\left\{
\begin{array}{lll}
\frac{\partial \Psi(x,t)}{\partial{t}} = \mathcal{A}\Psi(x,t), \quad \, t_{1}<t<t_{2},\\[2mm]
\Psi(t_{1}) = S(t_{1})\Psi_{0} + \Upsilon_{1}B.  \\
\end{array}
\right.
\end{equation*}
Since $ S(t_{1})\Psi_{0} + \Upsilon_{1}B \in \mathcal{D}(\mathcal{A}),$ Once again the unique classical solution of the previous equation is given by 
\begin{equation*}
\begin{array}{lll}
\Psi_{2}(t) = S(t)\Psi_{0}+S(t-t_{1})\Upsilon_{1}B,
\end{array}
\end{equation*}
furthermore, we have $\Psi_{2}(.)$ is left continuous at $t_{2}$ and $\Psi_{2}(t_{2}) \in D(\mathcal{A}).$ If we continue in the same manner then at the k-th step, in the interval $J_{k} = [t_{k-1}, t_{k}),$ $k \in \theta_{1}^{n+1}$ we get the following unique classical solution
$$\Psi_{k}(t) = S(t) \Psi_{0}+\displaystyle{\sum_{m=1}^{k-1}} S(t-t_{m})\Upsilon_{m}B .$$
Now, we define the function $\Psi$ as follows:
\begin{equation*}
\Psi(t) =
\left\{
\begin{array}{lll}
\Psi_{1}(t) & t_{0}< t <t_{1},&\\[2mm]
 & \vdots & \\
\Psi_{k}(t) & t_{k-1}< t <t_{k},&  \\
& \vdots & \\
 \Psi_{n+1} & t_{n}< t <T,&\\
\end{array}
\right.
\end{equation*}
\\
it is clear that $\Psi(.)$ is the unique classical solution of \eqref{s3}.
\end{proof} 
Note that if $\Psi_{0} \in \mathcal{H},$ the existence of a solution for the system \eqref{s3} can be proved in a similar way by assuming that $\Upsilon_{k} \in L^{2}(\Omega).$
\section{Impulse Controllability of the System \eqref{s3.1}}
Control theory deals with how an arbitrary initial state can be directed exactly or approximately close to a given final state using a set of admissible controls. In this section, we study the impulse approximate controllability for the system below with one pulse $\tau \in (0,T)$. Here, the control function acts on a subdomain $\omega$ and at one point of time $\tau \in (0,T)$ (see more in \cite{MBER}).\\
We consider the following wave equation with one pulse:\\
\begin{equation}\label{s3.1}
\left\{
\begin{array}{lll}
\frac{\partial^{2}\kappa(x,t)}{\partial{t}^{2}}-\frac{\partial^{2}\kappa(x,t)}{\partial{x}^{2}} = 0, &\text{ in }\,\,\,\, \Omega\, \times(0, T) \backslash\left\{\tau \right\},& \\
\kappa(x,\tau^{+}) = \kappa(x,\tau^{-}), &\text{ in }\,\,\,\, \Omega, &\\
\frac{\partial \kappa(x,\tau^{+})}{\partial{t}} -\frac{\partial \kappa(x,\tau^{-})}{\partial{t}} = \mathbf{1}_{\omega} \Upsilon_{\tau}(x), &\text{ in }\,\,\,\, \Omega,&\\
\kappa(0,t) = \kappa(1,t)= 0, &\text{ on }\,\,\,\, (0,T),&\\
\kappa(x,0) = \kappa^{0}(x) , \frac{\partial \kappa(x,0)}{\partial{t}}=\kappa^{1}(x),  &  \text{ in }\,\,\,\, \Omega.& \;\;\;\;\\
\end{array}
\right.
\end{equation}															
Let us first remark that \eqref{s3.1} is equivalent to the following system:

\begin{equation}\label{3.2}
\left\{\begin{array}{ll}
\frac{\partial^{2} \kappa(x, t)}{\partial t^{2}}-\frac{\partial^{2} \kappa(x, t)}{\partial x^{2}}=\mathbf{1}_{\omega} \Upsilon_{\tau}\left(x\right)\delta_{\tau}(t), & \text{ in } \Omega\times (0,T), \\
\kappa(0, t)=\kappa(1, t)=0, &\text{ on } (0,T), \\
\kappa(x, 0)=\kappa^{0}(x), \frac{\partial \kappa(x, 0)}{\partial t}=\kappa^{1}(x), & \text{ in }\Omega.
\end{array}\right.
\end{equation}
where $\delta_{\tau}$ is the Dirac measure at time $t = \tau$, which in turn is equivalent to the following abstract control problem

\begin{equation}\label{3.3}
\left\{
\begin{array}{lll}
\psi^{'}(t)= \mathcal{A} \psi(t) + B_{\omega} \Upsilon_{\tau}\delta_{\tau}(t), \quad &\text{ on }\,\,\,\, (0,T),&\\[2mm]
\psi(0) = \psi_{0},\\
\end{array}
\right.
\end{equation}
where
$\psi(0)=\left(\kappa^{0}, \kappa^{1}\right)^\top,\;
B_{\omega}=
\left(0,  \mathbf{1}_{\omega}\right)^\top$ and 
 $\mathcal{A} $  is an  unbounded linear operator defined on $\mathcal{H}$ as in \eqref{fi}.\\
Next, we introduce a lemma that will be frequently used in what follows: 
\begin{lem} \cite{Tu}
The  eigenfunctions of $\mathcal{A}$  are given by 
$$\Phi^{n}=\left(\sin(n\pi x), \lambda_n\sin(n\pi x)\right)^\top,\quad  \lambda_{n} = in\pi,\; n\in \mathbb{Z}^{*},$$
they forms an orthonormal basis in  $\mathcal{H}$. Moreover,  the group $(S(t))_{t\in \mathbb{R}}$ and its generator can be represented as follows:
\begin{eqnarray*}
 \mathcal{A}x &=& \sum_{n = 1}^{\infty} \lambda_{n} \langle x, \phi_{n} \rangle \phi_{n}, \; x \in \mathcal{H}\\
S(t)x &=& \sum_{n = 1}^{\infty} e^{- n^2 t}\langle x, \phi_{n} \rangle \phi_{n}, \; x \in \mathcal{H}.
\end{eqnarray*}
\end{lem}

\begin{defi}[Impulse approximate Controllability]
System \eqref{3.3} is said to be impulse approximately controllable on $[0,T],$ if for all desired state $\psi^{1} \in \mathcal{H},$ and initial condition $\psi_{0} \in \mathcal{H},$ and for  all $\varepsilon > 0 $ there exists $\Upsilon_{\tau} \in  L^2(\omega)$ such that the mild solution $\psi(t,\Upsilon_{\tau},\psi_0)$ of \eqref{3.3}  verifies:
\begin{equation*}
\left\|\psi(T,\Upsilon_{\tau},\psi_0) - \psi^{1}\right\|_{\mathcal{H}} <\varepsilon.
\end{equation*}
\end{defi}

The following lemma will be used to establish the main result of this section. We consider  $G : W \longrightarrow Z $ a linear bounded operator between Hilbert spaces  $W$ and $Z$.
\begin{lem}(see \cite{RFCAJP,RFCHJZ,HLNMJS})\label{exl}
The following statements are equivalent:
\begin{enumerate}
\item $ \quad \overline{\operatorname{Rang}\left(G\right)}=Z ,$\\[1.5mm]
\item $ \quad \ker\left(G^*\right) = \{0\}$\\[1.5mm]
\item  $\left\langle G G^* v , v\right\rangle>0,\quad v \neq 0 \quad\text{in}\,\, Z,$\\[1.5mm]
\item  $\lim _{\alpha \rightarrow 0^{+}} \alpha\left(\alpha I+G G^*\right)^{-1} v =0.$\,\, $\forall v \in Z.$
\end{enumerate}
\end{lem}
where $G^* $ is the adjoint operator of $G$. Now, we are in a position to state the main result of this section.
\begin{thm}
If $\omega = \Omega$, the system \eqref{s3.1} is approximate impulse controllable at any time $T > 0.$\\
If $\omega \varsubsetneq \Omega$  and $\{\Phi^1,\Phi^2,\ldots,\Phi^{N}\},\,N\in \N,$ forms  a finite basis of $ \mathcal{H}$, the system \eqref{s3.1} is approximate impulse controllable at any time $T > 0$. 
\end{thm}
\begin{proof}
According to Lemma \ref{exl}, see also  \cite[Theorem 2.43, p.56]{cor}, the approximate impulse controllability of the system \eqref{3.3} at $T>\tau,$ returns to the following uniqueness result for the  adjoint problem:
\begin{equation}\label{3.12}
\left\{\begin{array}{l}
\vartheta_{t}=-\mathcal{A} \vartheta \\
B_{\omega}^{*} \vartheta(x,\tau)=0 \quad \text{ a.e}\quad \text{in}\quad \Omega \quad 
\end{array}\right\} \Longrightarrow \vartheta \equiv 0.
\end{equation}
\\
Indeed, we shall distinguish two cases: $\omega = \Omega ;$ $\omega \varsubsetneq \Omega .$ \\
$\blacktriangleright$\textbf{Case 1:} For $\omega = \Omega ,$  
the solution of the adjoint problem with the following initial condition
$$\vartheta_{0}=\sum_{k \in \mathbb{Z}^{*}} a_{k} \lambda_{k} \Phi^{k} \in L^{2}(\Omega) \times H^{-1}(\Omega),$$
is given by
$$\vartheta(x,t)=\sum_{k \in \mathbb{Z^{*}}} a_{k} e^{-\lambda_{k} t} \lambda_{k} \Phi^{k}(x).$$
We assume that 
$$B_{\omega}^{*} \vartheta(x,\tau)=0, \quad \text{ a.e}\quad \text{in}\quad \Omega $$
that is
\begin{equation}\label{3.13}
 \sum_{k \in \mathbb{Z^{*}}} a_{k} e^{-\lambda_{k} \tau}\lambda_{k} \Phi^{k}(x) = 0,  \quad \text{a.e} \quad   \text{in}\quad \omega.
\end{equation}
then
$$ \left\|\sum_{k \in \mathbb{Z^{*}}} a_{k} e^{-\lambda_{k} \tau} \lambda_{k} \Phi^{k}(.)\right\|_{L^{2}(\Omega) \times H^{-1}(\Omega)} = 0, $$
since $(\lambda_{k} \Phi^{k})_{k\in \Z^{*}}$ is an orthonormal basis in  $L^{2}(\Omega) \times H^{-1}(\Omega)$,  we get

$$\sum_{k \in \mathbb{Z^{*}}} \left|a_k\right|^2  \left|e^{-\lambda_{k} \tau}\right|^{2}  = 0 ,$$
thus
$$\sum_{k \in \mathbb{Z^{*}}} \left|a_k\right|^2 = 0  ,$$
therefore, we obtain $a_{k} = 0$,   for all $k \in \mathbb{Z^{*}},$
which prove that
$$ \vartheta \equiv 0.$$

$\blacktriangleright$\textbf{Case 2:} For $\omega \varsubsetneq \Omega ,$  and $\{\Phi^1,\Phi^2,\ldots,\Phi^{N}\}$, form a finite basis of $H$, the solution of the adjoint problem with the following initial  condition:
$$\vartheta_{0}=\sum_{k =1}^{N} a_{k}\lambda_{k} \Phi^{k} ,$$
is given by
$$\vartheta(x,t)=\sum_{k =1}^{N} a_{k} e^{-\lambda_{k} t}\lambda_{k} \Phi^{k}(x).$$
We assume that
\begin{equation*}
\displaystyle B_{\omega}^{*}\,\vartheta(x,\tau)=0 \quad \text{a.e} \quad \Omega,
\end{equation*}
that is 
\begin{equation}\label{3.4} 
\sum_{k=1}^{N} a_{k}\lambda_{k}^2\, e^{-\lambda_{k} \tau} \sin(k\pi x) =0 \quad \text{a.e} \quad \omega. 
\end{equation}
We recall that the classic Chebychev polynomials of second species are given by the following  relation:
\begin{equation*}
\left\{\begin{array}{ll}
U_0=1,\;U_1=2X,\\
\\
U_{n+1}=2 X U_{n}-U_{n-1}, \quad \forall n \geq 1 .\\
\end{array}\right.
	\end{equation*}
Polynomials $(U_{n})$ can be defined alternatively by  trigonometric forms of their associated polynomial functions on $\omega.$ Indeed, for every $n\in \N^{*}$ we have 
$$\sin(n\pi x) = \sin(\pi x) U_{n-1}(\cos(\pi x)), \quad \forall x \in \omega.$$
Therefore the assumption $\eqref{3.4}$, becomes
\begin{equation}\label{3.5} 
\sum_{k=1}^{N} a_{k}\lambda_{k}^2\, e^{-\lambda_{k} \tau}  \sin(\pi x) U_{k-1}(\cos(\pi x)) =0, \quad \text{a.e} \quad \omega.
\end{equation}
Hence, it is enough to prove the following implication:
\begin{equation}\label{3.16}
\sum_{k=1}^{N} a_{k}\lambda_{k}^2 e^{-\lambda_{k} \tau}  \sin(\pi x) U_{k-1}(\cos(\pi x)) =0 \quad \text{a.e} \quad \omega \Longrightarrow a_{k} = 0 \quad k=1,2,...,N .
\end{equation}
For $N=1$  the implication \eqref{3.16} is clearly verified. Now, we assume that it is true for $n < N,$ and we prove it for $n+1$. By taking
\begin{equation}\label{3.6} 
\sum_{k=1}^{n+1} a_{k}\lambda_{k}^2 e^{-\lambda_{k} \tau}  \sin(\pi x) U_{k-1}(\cos(\pi x)) =0 \quad \text{a.e} \quad \omega.
\end{equation}
since $\sin(\pi x) \neq 0$, for all $x \in \omega$, we obtain that 
\begin{equation}\label{tsou}
a_{1}\lambda_{1}^2 e^{-\lambda_{1} \tau}+...+a_{n+1}\lambda_{n+1}^2 e^{-\lambda_{n+1} \tau} U_{n}(\cos(\pi x))=0, \quad \text{a.e} \quad \omega. 
\end{equation}
we compute the derivative of (\ref{tsou}), we obtain that
\begin{equation}\label{3.19}
 a_{2}\lambda_{2}^2 e^{-\lambda_{2} \tau} \pi U_{1}^{(1)}(\cos(\pi x))+\ldots+a_{n+1}\lambda_{n+1}^2 e^{-\lambda_{n+1} \tau} \pi U_{n}^{(1)}(\cos(\pi x))=0, \quad \text{a.e} \quad \omega, 
\end{equation}
again, if we compute the derivative of the previous formula \eqref{3.19}, we obtain that
\begin{equation*}
a_{3}\lambda_{3}^2 e^{-\lambda_{3} \tau} \pi^{2} U_{2}^{(2)}(\cos(\pi x))+...+a_{n+1}\lambda_{n+1}^2 e^{-\lambda_{n+1} \tau} \pi^{2} U_{n}^{(1)}(\cos(\pi x))=0 \quad \text{a.e} \quad \omega,
\end{equation*}
we continue in the same manner. Finally, at the n-$\mathrm{th}$ step, we get
$$a_{n+1}\lambda_{n+1}^2 e^{-\lambda_{n+1} \tau} \pi^{n} U_{n}^{(n)}(\cos(\pi x))=0 \quad \text{a.e} \quad \omega, $$
since  $U_{n}^{(n)}(\cos(\pi x)) \in \mathbb{R}^{*},$ then $a_{n+1} = 0.$\\
This ends to the desired implication.
\end{proof}

\section{Impulse Observability Inequality for \eqref{s3.1}}
In this section we make use of the strategy presented by Lions \cite{lions} to obtain an observation estimate at one instant of time $\tau=2,$\,( $0< \tau  <T)$ for the impulsive wave equation \eqref{s3.1}. \\
For $(\Phi^{0},\Phi^{1})\in C^{\infty}_{0}(\Omega)\times C^{\infty}_{0}(\Omega),$ the homogeneous problem associated to \eqref{s3.1} is  as follows:
\begin{equation}\label{4.1}
\left\{
\begin{array}{lll}
\frac{\partial^{2}\Phi(x,t)}{\partial{t}^{2}}-\frac{\partial^{2}\Phi(x,t)}{\partial{x}^{2}} = 0, &\text{in}\,\,\,\, \Omega\, \times(0, T), \\
\Phi(0,t) = \Phi(1,t)= 0, &\text{on}\,\,\, (0,T),\\
\Phi(x,0) = \Phi^{0}(x) , \frac{\partial \Phi(x,0)}{\partial{t}}=\Phi^{1}(x),  & \text{in}\,\,\, \Omega, \;\;\;\;
\end{array}
\right.
\end{equation}
System \eqref{4.1} admits a unique solution expressed in Fourier series as
$$ \Phi(x,t)=\sum_{n=1}^{\infty}\left(a_{n} \cos (n \pi t)+\frac{b_{n}}{n \pi} \sin (n \pi t)\right) \sin (n \pi x) ,$$
where
\begin{equation*}
\Phi^{0}(x) = \sum_{n=1}^{\infty} a_{n}\sin( n \pi x)\quad\text{and}\quad
\Phi^{1}(x) = \sum_{n=1}^{\infty} b_{n}\sin( n \pi x),
\end{equation*}
$(a_{n})$ and $(b_{n})$ are the coefficients of Fourier in the orthogonal basis of $L^{2}(\Omega)$
$$ \theta_{n}(x) = \sin(n\pi x),\quad n=1,2,\ldots.$$
Next, we consider the following backward problem
\begin{equation}
\left\{
\begin{array}{lll}\label{4.2}
 \frac{\partial^{2} \Psi(x,t)}{\partial{t}^{2}} -\frac{\partial^{2} \Psi(x,t)}{\partial{x}^{2}} = 0, &\text{ in }\,\,\,\, \Omega\, \times(0, T) \backslash\left\{\tau\right\},\\
\Psi(x,\tau^{-}) = \Psi(x,\tau^{+}), &\text{ in }\,\,\,\, \Omega, &\\
\frac{\partial \Psi(x,\tau^{+})}{\partial{t}} -\frac{\partial \Psi(x,\tau^{-})}{\partial{t}} = -\mathbf{1}_{\omega}\Phi_{t}(x,\tau)\delta_{\tau}, &\text{ in }\,\,\,\, \Omega ,\\
\Psi(0,t)= \Psi(1,t)= 0,  &\text{ on }\,\,\,\, (0,T),\\
\Psi(T) = \Psi^{'}(T) =0, & \text{ in }\,\,\,\, \Omega. \;\;\;\;\\
\end{array}
\right.
\end{equation}
We can easily check that this system is well-posed, by a simple calculation we obtain
$$\int\int_{Q} \frac{\partial^{2} \Psi(x,t)}{\partial{t}^{2}}\Phi(x,t) dxdt-\int\int_{Q}\frac{\partial^{2} \Psi(x,t)}{\partial{x}^{2}}\Phi(x,t) dxdt = - \int\int_{Q}\mathbf{1}_{\omega} \Phi_{t}^{2}(x,\tau)dxdt, $$
with $Q= \Omega\times (0.T).$\\ 
A double integration by parts yields
$$T \int_{\omega}\Phi_{t}^{2}(x,\tau)dx  =  \int_{\Omega} \Phi^{0}\Psi^{'}(0)-\Psi(0)\Phi^{1}dx ,$$
we define the operator $\Lambda$ as follows:
 $$\Lambda(\Phi^{0},\Phi^{1}) = (\Psi^{'}(0),-\Psi(0)) \quad \forall  (\Phi^{0},\Phi^{1})\in C^{\infty}_{0}(\Omega)\times C^{\infty}_{0}(\Omega),$$
then, $$ T \int_{\omega}( \Phi_{t}(x,\tau)^{2} dx = <\Lambda(\Phi^{0},\Phi^{1}),(\Phi^{0},\Phi^{1})>_{L^{2}\times L^{2}} .$$
For  $(\Phi^{0},\Phi^{1}) \in  C^{\infty}_{0}(\Omega)\times C^{\infty}_{0}(\Omega),$ we define the semi-norm 
$$\left\| \lbrace\Phi^{0},\Phi^{1} \rbrace\right\|_{F} := T^{\frac{1}{2}}\left(\int_{\omega}\mid\Phi_{t}(x,\tau)\mid^{2}dx\right)^{\frac{1}{2}}. $$
In what follows, we will prove the impulse observability inequality (see \cite{Tu,ABGDMCDSM,ETGWYX}), which consist of the existence of a constant $c>0$ such that for all $\left( \Phi^{0},\Phi^{1} \right) \in L^{2}(\Omega)\times H^{-1}(\Omega)$ the mild solution $\Phi$ of the problem \eqref{4.1} satisfies
\begin{equation}\label{4.3}
\int_{\omega}\mid \Phi_{t}(x,\tau)\mid^{2}dx \geq c\, \| \lbrace\Phi^{0},\Phi^{1} \rbrace\|^{2}_{ L^{2}(\Omega)\times H^{-1}(\Omega)} .
\end{equation}
Indeed, we shall distinguish two cases: $\omega = \Omega ;$ $\omega \varsubsetneq \Omega .$ \\
$\blacktriangleright$\textbf{case 1:} For $\omega = \Omega,$
we may express $\displaystyle\int_{\Omega}\mid \Phi_{t}(x,\tau)\mid^{2}dx $  in terms of the Fourier coefficients $(a_{n})$ and $(b_{n})$ as follows:
\begin{equation*}
\begin{array}{lll}
\displaystyle
\int_{\Omega}\mid\Phi_{t}(x,\tau)\mid^{2}dx &=&\displaystyle\frac{1}{2} \sum_{n=1}^{\infty} \left( n\pi \,a_{n}\cos(n\pi \tau)+ b_{n}\sin(n\pi \tau) \right)^{2}\\[5mm]
&=&\displaystyle\frac{1}{2}\sum_{n=1}^{\infty}\left((n\pi a_{n})^{2} +b_{n}^{2}\right)\sin^{2}(n\pi \tau +y_{n}),
\end{array}
\end{equation*}
with  $y_{n}$ satisfied  $$\sin(y_{n})= \frac{b_{n}}{\sqrt{(n\pi a_{n})^{2} +b_{n}^{2}}} \quad \text{and} \quad \cos(y_{n})= \frac{n\pi a_{n}}{\sqrt{(n\pi a_{n})^{2} +b_{n}^{2}}}.$$
On the other hand
\begin{equation*}
\begin{array}{lll}
\displaystyle
\left\|\phi_{0}\right\|_{L^{2}(\Omega)}^{2}+\left\|\phi_{1}\right\|_{H^{-1}(\Omega)}^{2} = \sum_{n=1}^{\infty}\left(a_{n}^{2}+\frac{b_{n}^{2}}{n^{2}
 \pi^{2}} \right).
\end{array}
\end{equation*}
So, it's enough to prove the existence of a positive constant $c$ such that
$$\displaystyle\sum_{n=1}^{\infty}\left((n\pi a_{n})^{2} +b_{n}^{2}\right)\sin^{2}(n\pi \tau +y_{n})\geq c \displaystyle  \sum_{n=1}^{\infty}\left(a_{n}^{2}+\frac{b_{n}^{2}}{n^{2}
 \pi^{2}} \right),$$
this inequality would require a lower bound of the form
\begin{equation}\label{4.23}
\begin{array}{lll}
\displaystyle
\vert\sin(n\pi \tau +y_{n})\vert \geq c \quad \forall n \in \N^{*},
\end{array}
\end{equation}
since $\tau = 2,$ the inequality \eqref{4.23} is equivalent to the following:
\begin{equation}\label{4.25}
\begin{array}{lll}
\displaystyle
\vert\sin(y_{n})\vert \geq c \quad \forall n \in \N^{*},
\end{array}
\end{equation}
this inequality is false for all $y_{n}.$ Indeed, if $y_{n},$ is
 expressed as
$$y_{n} = n \pi \quad n \in \N^{*}, $$
then, foll $n$ peer
$$\sin(y_{n}) = \sin( n \pi)= 0 .$$
In this case, $\sin(y_{n})= 0,$ for an infinite number of values of $n.$ Thus, inequality \eqref{4.25} cannot be true. But, for certain values of $y_{n}$ this inequalities may be obtained. For instance, if $  b_{n} = a_{n}$ for all $n \in \N,$ 
\begin{equation*}
\displaystyle
\vert\sin(y_{n})\vert = \frac{\vert b_{n}\vert}{n\pi\sqrt{b_{n}^{2} +\frac{b_{n}^{2}}{(n\pi)^{2}}}}, 
\end{equation*}
therefore 
\begin{equation*}
\vert\sin(y_{n})\vert \geq \frac{1}{n\pi \sqrt{2}},
\end{equation*}
and this is the best lower bound one may expect, this implies that
\begin{equation*}
\begin{array}{lll}
\displaystyle
\displaystyle\sum_{n>0}^{\infty}\left((n\pi a_{n})^{2} +b_{n}^{2}\right)\sin^{2}(n\pi \tau +y_{n})\geq \frac{1}{2} \displaystyle\sum_{n>0}^{\infty}\left( a_{n}^{2} + \frac{b_{n}^{2}}{n^{2}\pi^{2}} \right),
\end{array}
\end{equation*}
which is equivalent to 
\begin{equation*}
\int_{\Omega}\mid \Phi_{t}(x,\tau )\mid^{2}dx \geq \frac{1}{2}\, \| \lbrace\Phi^{0},\Phi^{1} \rbrace\|^{2}_{ L^{2}(\Omega)\times H^{-1}(\Omega)} .
\end{equation*}

$\blacktriangleright$\textbf{case 2:} For $\omega \varsubsetneq \Omega ,$ in this case, we'll do a numerical check  for the class of initial conditions $ a_n = b_n ,$ for which we had proved the observability inequality \eqref{4.3} in case $\omega = \Omega$.  In other words, let's take $\omega = ]0,\frac{1}{2}[ \varsubsetneq \Omega,$ and we put the following initial conditions:
\begin{equation*}
\Phi^{0}_{N}(x) = \sum_{n=1}^{N} a_n \sin( n \pi x)\quad\text{and}\quad
\Phi^{1}_{N}(x) = \sum_{n=1}^{N} b_n \sin( n \pi x),
\end{equation*}
we note  $\Phi^{N}$ the solution of \eqref{4.1} corresponding to the above initial conditions.
The figures below shows the variation of the quantity  $\eqref{4.5}$ vis-a-vis $N,$ for the following particular initiales conditions: $a_n = b_n = k \in \mathbb{R}; $   $a_n=  b_n = n;$ and $a_n = b_n =  n \pi.$
\begin{equation}\label{4.5}
\frac{\displaystyle\int_{\omega}\mid \Phi^{N}_{t}(x,t_{1})\mid^{2}dx}{\displaystyle \| \lbrace\Phi^{0}_{N},\Phi^{1}_{N} \rbrace\|^{2}_{ L^{2}\times H^{-1}}} .
\end{equation}

\begin{tabular}{cc}
\includegraphics[scale=0.6]{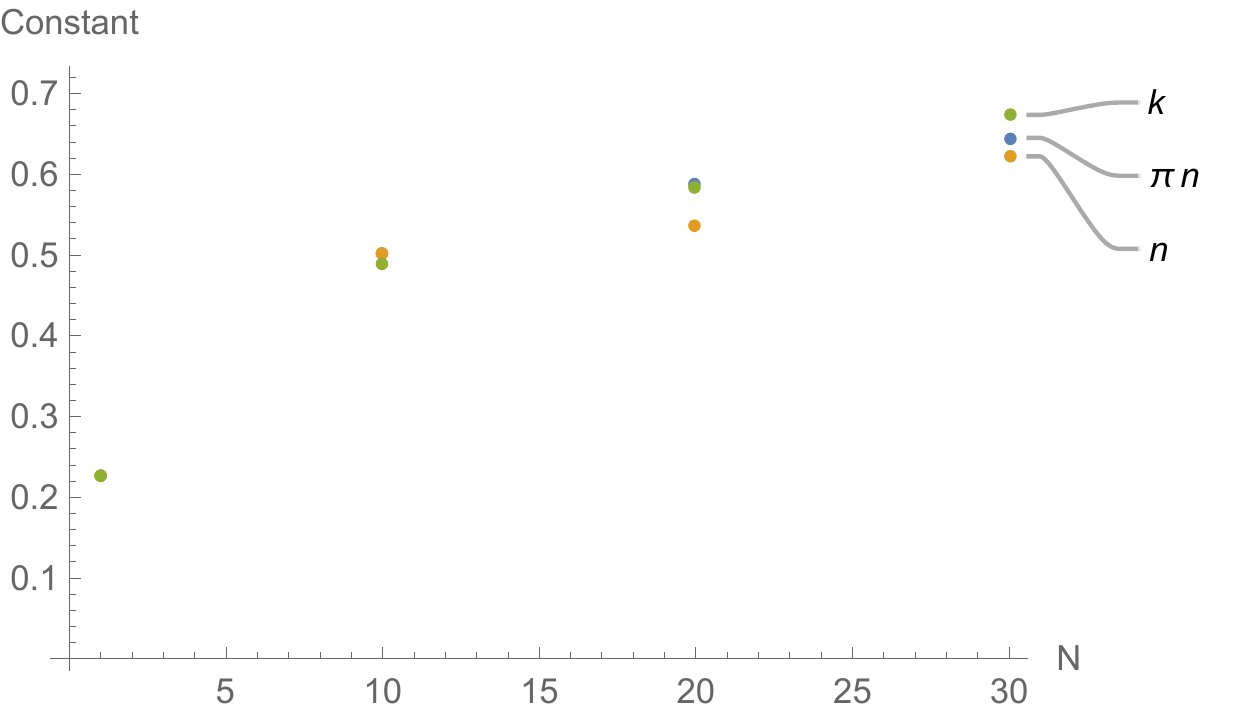} &
\includegraphics[scale=0.6]{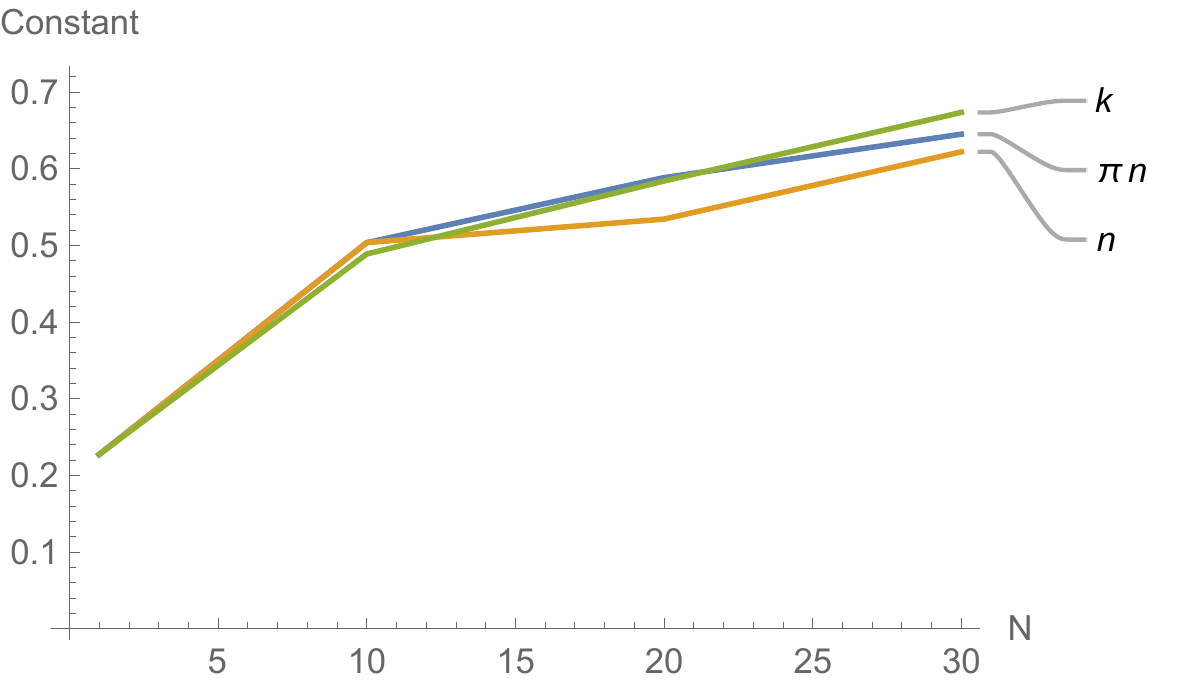} \\
\end{tabular}
In these figures, we notice that \eqref{4.5} is always bigger than $c \simeq 0,227$ and it's strictly increasing, which validates the inequality of observability  \eqref{4.3} for these particular cases.


\section{Conclusion and Open problems}
Impulse control is a very weak control function, it acts only at one instant of time, and the problem of the wave equation with impulse control on a subdomain $\omega$, without any geometrical conditions on the space $\Omega$ is very challenging. In this paper, we made an analysis of this problem from different angles and we obtained some interesting results. But there are still several open questions that we have not been able to answer like the problem of controllability of the system \eqref{s3.1}, in a case $\omega \varsubsetneq \Omega ,$ in infinite dimension. The novelty in this work is that the literature on impulse control systems is very small, there are very few numbers of works on systems with impulse control, to our knowledge wave equation with impulse control have not been studied yet.


\begin{thebibliography}{10}
\bibitem{BMJ} 
\newblock M. Benchohra, J. Henderson, and S. Ntouyas,
\newblock Impulsive differential equations and inclusions,
\newblock \emph{New York: Hindawi Publishing Corporation,} \textbf{2} 2006.

\bibitem{ABGDMCDSM}
\newblock A. Bensoussan, G. D. Prato,  M.C. Delfour and S. Mitter,
\newblock Representation and Control of Infinite Dimensional Systems,
\newblock \emph{Birkh{\"a}user Basel,} 2007.

\bibitem{CDGHL} O. Camacho and H. Leiva,
\newblock Impulsive semilinear heat equation with delay in control and in state,
\newblock\emph{ Asian Journal of Control,} \textbf{22} (2020), 1075-1089.

\bibitem{ACCDGHL} 
\newblock A. Carrasco, G. Guevara and H. Leiva,
\newblock Controllability of the impulsive semilinear beam equation with memory and delay,
\newblock\emph{IMA Journal of Mathematical Control and Information,} \textbf{36} (2019),  213-223.

\bibitem{cor}
\newblock J. M. Coron,
\newblock Control and nonlinearity,
\newblock \emph { American Mathematical Society, Boston,} 2007.

\bibitem{RFCAJP}
\newblock R. F. Curtain and A. J. Pritchard,
\newblock Infinite Dimensional Linear Systems Theory,
\newblock \emph{Springer-Verlag Berlin Heidelberg,}  1978.

\bibitem{RFCHJZ} 
\newblock R. F. Curtain and H. Zwart,
\newblock An introduction to infinite-dimensional linear systems theory, Vol.\textbf{21}
\newblock \emph{Springer-Verlag, New York,} 2012.

\bibitem{CDJUHLOC} 
\newblock C. Duquea, J. Uzcategui, H. Leiva and O. Camacho
\newblock Approximate controllability of semilinear strongly damped wave equation with impulses, delays, and nonlocal condi-tions
\newblock \emph{Journal of Mathematics and Computer Science,} \textbf{20} (2019), 108-121.

\bibitem{GL}
\newblock C. Guevara and H. Leiva,
\newblock Controllability of the impulsive semilinear heat equation with memory and delay,
\newblock\emph{Journal of Dynamical and Control Systems,} \textbf{24} (2018),  1-11.


\bibitem{KAY}
\newblock A.Y. Khapalov,
\newblock Exact controllability of second-order hyperbolic equations with impulse controls,
\newblock\emph{ Applicable Analysis,} \textbf{63} (1996), 223-238.

\bibitem{LHWZME}
\newblock H. Leiva, W. Zouhair and M. e. Entekhabi,
\newblock Approximate Controllability of semi-linear Heat equation with Non-instantaneous Impulses, Memory and Delay,
\newblock \emph{arXiv preprint,} (2020), \href{arXiv:2008.02094}{arXiv:2008.02094 }.

\bibitem{HLNMJS}
\newblock H. Leiva, N. Merentes, and J. Sanchez,
\newblock A characterization of semilinear dense range operators and applications,
\newblock \emph{Abstract and Applied Analysis,} \textbf{2013} (2013), 1-11.

\bibitem{lions} 
\newblock J.L. Lions,
\newblock Contr\^{o}labilit\'{e} exacte, perturbations et stabilisation de syst\`{e}mes distribu\'{e}s. Vol. \textbf{1} ,
\newblock \emph{ Masson, Paris,} 1988.

\bibitem{liu}
\newblock X. Liu,
\newblock Practical stabilization of control systems with impulse effects,
\newblock\emph{Journal of mathematical analysis and applications,} \textbf{166} (1992) 563-576.

\bibitem{farzana} 
\newblock F.A. McRae,
\newblock Practical stability of impulsive control systems,
\newblock \emph{Journal of mathematical analysis and applications,} \textbf{181} (1994),656-672.




\bibitem{MVDM}
\newblock V. D. Milman and A. D. Myshkis,
\newblock On the stability of motion in the presence of impulses,
\newblock \emph{Sibirskii Matematicheskii Zhurnal,} \textbf{1}(1960), 233-237.


\bibitem{MBER}
\newblock B. M. Miller and E. Y. Rubinovich,
\newblock Impulsive control in continuous and discrete-continuous systems,
\newblock \emph{Springer US, New York,} 2003.

\bibitem{pazy} 
\newblock A. Pazy,
\newblock Semigroups of linear operators and applications to partial differential equations. \newblock \emph{Springer Science and Business Media,} \textbf{44} 2012.

\bibitem{KDP} K.D. Phung, G. Wang, and Y. Xu,
\newblock Impulse output rapid stabilization for heat equations,
\newblock\emph{Journal of Differential Equations,} \textbf{263} (2017), 5012-5041.

\bibitem{SQG}
\newblock S. Qin and G. Wang,
\newblock Controllability of impulse controlled systems of heat equations coupled by constant matrices,
\newblock\emph{Journal of Differential Equations,} \textbf{263} (2017), 6456-6493.


\bibitem{yao}
\newblock A.M. Samoilenko and N.A. Perestyuk,
\newblock Impulsive differential equations,
\newblock\emph{world scientific,singapore,} 1995.



\bibitem{Tu}
\newblock M. Tucsnak and G. Weiss, 
\newblock Observation and control for
 operator semigroups,
\newblock\emph{Birkh{\"a}user, Basler,} 2009.


\bibitem{ETGWYX}
\newblock E. Tr{\'e}lat, G. Wang and  Y. Xu,
\newblock Characterization by observability inequalities of controllability and stabilization properties,
\newblock \emph{Pure and Applied Analysis,} \textbf{2} (2020), 93-122.

\end{thebibliography}
\end{document}